\documentclass[
  journal=small,
  manuscript=article-type,  
  year=2020,
  volume=37,
]{cup-journal}

\usepackage{amsmath}
\usepackage{amssymb}
\usepackage[nopatch]{microtype}
\usepackage{booktabs}
\usepackage[english]{babel}
\usepackage{csquotes}
\newtheorem{theorem}{Theorem}[section]

\title{On the equidistribution properties of patterns in prime numbers Jumping Champions, meta analysis of properties as Low-Discrepancy Sequences, and some conjectures based on Ramanujan's master theorem and the zeros of Riemann's zeta function}

\author{A. Ortiz-Tapia}
\affiliation{Universidad Abierta y a Distancia de M\'exico}
\email[F. Author]{egl.arturo.ortizta@unadmexico.mx} %


\begin{document}

\begin{abstract}
The Paul Erd\H{o}s-Tur\'an inequality is used as a quantitative form of Weyl' s criterion, together with other criteria to asses equidistribution properties on some patterns of sequences that arise from indexation of prime numbers, Jumping Champions (called here and in previous work, "meta-distances" or even md, for short). A statistical meta-analysis is also made of previous research concerning meta-distances to review the conclusion that meta-distances can be called Low-discrepancy sequences (LDS), and thus exhibiting another numerical evidence that md's are an equidistributed sequence. Ramanujan's master theorem is used to conjecture that the types of integrands where md's can be used more succesfully for quadratures are product-related, as opposite to addition-related. Finally, it is conjectured that the equidistribution of md's may be connected to the know equidistribution of zeros of Riemann's zeta function, and yet still have enough "information" for quasi-random integration ("right" amount of entropy).    \end{abstract}


\section{Introduction}
\subsection{Jumping Champions and Metadistances}
In their paper [9] of 1999 A. Odlyzko, M. Rubinstein, and M. Wolf made the following assertion

Wolf, Odlyzko, Rubinstein Conjecture 1.2
Conjecture 1.2. The jumping champions tend to infinity. Furthermore, any fixed prime p divides all
sufficiently large jumping champions.\cite{odlyzko1999jumping}

 the first assertion of Conjecture 1.2 was proved in 1980 by Erd˝os
and Straus \cite{erdos1980remarks},under the assumption of the truth of the Hardy-Littlewood prime pair conjecture \cite{hardy1923some}.

[...] Thus, 6 may not be the Jumping champion as x tends to infinity, but according to Conjecture 1.2, it can always be possible to find 6 as a difference between two prime numbers.


Among the efforts to study the patterns within prime numbers, it is to describe the ``Jumping Champions'' \cite{odlyzko1999jumping} which are the most frequent differences between prime numbers, that is to say
\begin{equation}\label{jumpers} 
\mathbb{J}={\max}[P_{i+1}-P_i]
\end{equation}
With the exception of the distance between the second and the first primes, namely, 2 and
3, all the other distances between prime numbers are even \cite{guy2004unsolved}.

Because 6 is the jumping champion up to about $n \approx 1.74 \times 10^35$ \cite{odlyzko1999jumping}, it was decided to investigate if the number 6 was associated systematically with other distances for subsets of this last range \cite{ortiz2008some}, and the exploration was performed for dyads, triads and and other $n$-ads. Because the dyads could be found more often \cite{ortiz2008some}, it was decided
to concentrate on them for further analysis. Out of the first one million distances between prime numbers, the dyad \{4, 6\} was found 15,860 times and the dyad \{6, 6\} was found 17,546 times, with diverse spacings between each dyad \cite{ortiz2008some}. Specifically, let the dyad
$\{6, 6\}$ be formed by:
\begin{equation}
\{P{i+1}-P_i,\;  P_{i+2} - P{i+1}\}   
\end{equation}

where this set of `\,`\emph{P}'s'\,' represent those prime numbers
capable of generating the required dyad. The set of distances between any two dyads like $\{6, 6\}$ was found with the following observation of the indices of primes related to those dyads:
\begin{eqnarray}
    \overbrace{i,i+2}^{j},\overbrace{i+k,i+k+2}^{j+1},\; k\in \mathbf{N}\\\nonumber
    md=(i+k_{j+1})-(i+k_{j+2}), \; md \in \mathbf{Md}\\\nonumber
\end{eqnarray}
where $md$ stands as shorthand for ‘meta-distances’, the way those distances were called in previous work \cite{ortiz2008some}, implying distances embedded in the set of distances between prime numbers, and alsoe $\mathbf{Md}$ is the set of all the meta-distances within a given range of prime numbers. There was a brief exploration to find a pattern within the set $\mathbf{Md}$, and none was found. Precisely it was the lack of apparent pattern in the meta-distances which lead to an inquiry whether they form a uniformly distributed sequence, which could be a low-discrepancy sequence, and in consequence the meta-distances could be used for multidimensional integration \cite{ortiz2008some}.

Several numerical tests were made using $\mathbf{Md}$, and also other sequences, all of them compared against the Gauss-Kronrod method implemented in Mathematica (GK) \cite{wolfram1999mathematica} or the analytical solution. The other sequences were created from the following generators:
\begin{itemize}
    \item Mathematica’s® Monte Carlo method of integration (MMC).
    \item A Monte Carlo method, using the random number generator of Mathematica (RMN) \cite{wolfram1999mathematica}, which in turn uses the Marsaglia-Zaman subtract-with-borrow generator for real numbers \cite{marsaglia1990toward}
    \item Prime intervals or first differences of primes (PI).
    \item The digits of the number $\pi$ (PiDi).
    \item The first $20,000$ digits of the number $e = 2.718281\cdots$ (eDi)
    \item The digits of $\sqrt{71}$ (71Di).
\end{itemize}
All the numerical bases based on the digits of some number were normalized in a similar fashion to the prime meta-distances: because the maximal number in those bases is the digit 9, all of them were divided by this digit and multiplied by two. Integration was done also by a convolution \cite{ortiz2008some}, and was carried for different types of integrands, with different types of bounds and integrating in one or more dimensions, up to four dimensional integrands. The results were analyzed statistically using the mean and standard deviation of the relative error, whereupon $\mathbf{Md}$ resulted in $9\% \leq \bar{x} \leq 11\%$ and $2\%\leq \sigma \leq 25\%$. The integration using the other types of sequences resulted in similar variability \cite{ortiz2008some}.

As for the speed of convergence, almost all sequences had at least linear convergence, that is to say, suppose we have a sequence \(\{x_n\}\) such that \(x_n\rightarrow x_\infty\) in \(\mathbb{R}^k\). We say the convergence is \emph{linear} if there exists \(r\in(0, 1)\)
\begin{equation}
\frac{||x_{n+1}-x_\infty||}{||x_n-x_\infty||}\leq r    
\end{equation}

for all \(n\) sufficiently large.
In fact, $\mathbf{Md}$ consistently converged with at most $15,000$ evaluations of the integrand, regardless of the dimension, as opposite to all other sequences, converging with at least $20,000$ integrand evaluations \cite{ortiz2008some}.

\subsection{Leveque, Wozniakowski, and Halton theorems}

In previous work \cite{ortiz2008probability} there was some exploration of Leveque's inequality \cite{kuipers2012uniform}, Wozniakoski's theorem \cite{wozniakowski1991average}, and Halton's theorem \cite{halton1960efficiency}, concerning the identification of a discrepancy sequence as such. Leveque's inequality resulted in an average of $D_N\approx 0.8$, Wozniakoski's rendered $T_N^* \approx 0.2$ and Halton theorem was an average of $\approx9$, accross dimensions of integration $s\in {1,2,3,4,5,6,7}$. No comparison was made with other sequences. Both Leveque's inequality and Wozniakoski theorem agree in order of magnitude about the discrepancy of $\mathbf{Md}$, whereas Halton theorem would agree in order of magnitude with the average error across all dimensions, which was $\bar{\varepsilon} \approx 4.9$ \cite{ortiz2008probability}. Again, no comparison was made with other sequences.

In this work new comparisons are made, and new evidence is presented for considering $\mathbf{Md}$ as an Equidistributed sequence.

\section{Ramanujan's master theorem and Conjecture of types of integrand better suited for $\mathbf{MD}$}

In \cite{BanffMultiplicativeDistributions} it was shown that kernels of "multiplicative" type of Fredholm integral equations are better suited for being solved with Low-Discrepancy Sequences, specifically with $\mathbf{Md}$. Multiplicative, in the sense that the kernel is composed mostly of products, as opposite of having sums inside the kernel. The reason for this numerical behavior was conjectured to be related to the higher probability of divisibility of numbers as opposite to the probability of addition. Although this conjecture was substantiated with some numerical exploration, a stronger argument can be found using Ramanujan' s Master Theorem, named after Srinivasa Ramanujan \cite{amdeberhan2012ramanujan}.

Ramanujan' s Master Theorem is a technique that provides an analytic expression for the Mellin transform of an analytic function. The result is stated as follows : 
  \begin{theorem}[Ramanujan's Master theorem]
\label{Theorem:RamanujanMasterTheorem2}
  If a complex-valued function $f (x)$ has an expansion of the form
  \begin{equation}
      f (x) = \sum _ {k =      0}^{\infty} {\frac {\, \varphi (k)\,} {k!}} (-x)^{k}
  \end{equation}     
  then the Mellin transform of $f (x)$ is given by
\begin{equation}
    \int _ {0}^{\infty} x^{s - 1} f (x)\, dx = \Gamma (s)\, \varphi (-s)
\end{equation}\label{Eq:MellinTransformGammaF}     
  \end{theorem}
 
Since the Gamma function (the generalized factorial) is present in this theorem, any analytic function that complies with Theorem \ref{Theorem:RamanujanMasterTheorem2} would be expressed as some type of product, thus it can be conjectured that $\mathbf{Md}$ solve both integrals or integral equations with integrands or kernels respectively that comply with theorem \ref{Theorem:RamanujanMasterTheorem2}

\pagebreak

\section{Statistics}
The rule of thumb seems to be: If the skewness is between -0.5 and 0.5, the data are fairly symmetrical.	If the skewness is between -1 and – 0.5 or between 0.5 and 1, the data are moderately skewed. If the skewness is less than -1 or greater than 1, the data are highly skewed \cite{Skewness}. Thus, the results in Table \ref{Table:errors} would suggest that $\mathbf{Md}$ has the most symmetrical distribution of errors. 

It is also known that Positive kurtosis indicates that the data exhibit more extreme outliers than a normal distribution, negative kurtosis indicates that the data exhibit less extreme outliers than a normal distribution, and for a normal distribution, the value of the kurtosis statistic is zero \cite{Kurtosis}. Thus, under this interpretation of kurtosis, the results shown in table \ref{Table:errors} would suggest that $\mathbf{Md}$ would exhibit the least amount of extreme outliers.

\begin{figure}[h]
\centering\includegraphics[width=10cm]{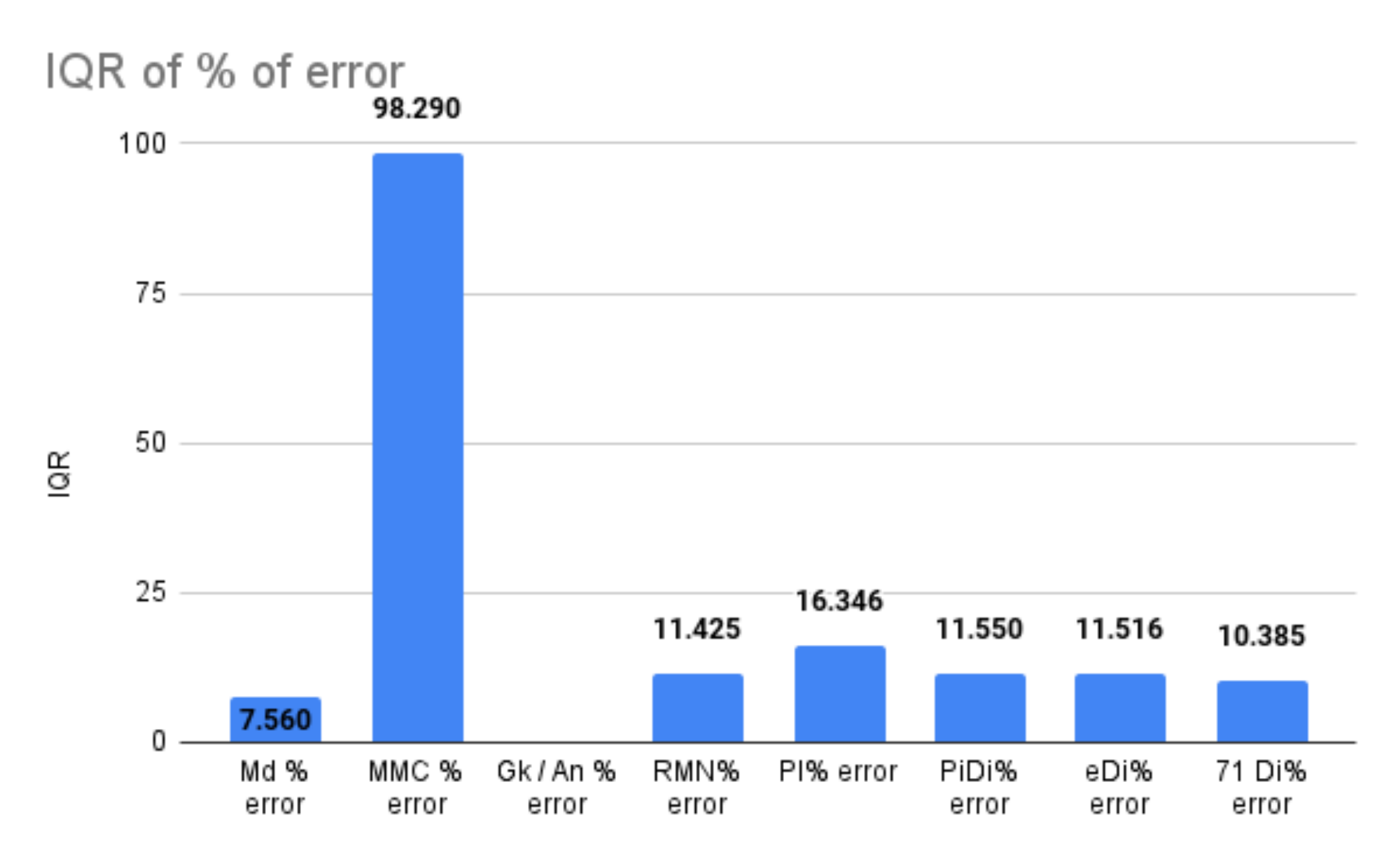}
\caption{Comparison of Interquartile Range of the percentage of errors for the different sequences}\label{Figure:IQR}
\end{figure}

As for the Interquartile Range, the comparison shown in Fig.\ref{Figure:IQR} shows evidence that $\mathbf{Md}$ has the least dispersion, albeit the standard deviation of its errors shows the opposite in Table \ref{Table:errors}.

\begin{table}[hbt!]
\begin{threeparttable}
\caption{\% error (excluding outliers)}
\label{Table:errors}
\begin{tabular}{lrrrrrrrrr}\toprule
&\textbf{Md} &MMC  &Gk / An  &RMN &PI &PiDi &eDi &71 Di \\\midrule
average &7.6 &1.0 &NA &6.3 &14.4 &7.9 &7.2 &7.1 \\
stdv &4.7 &1.0 &NA &8.1 &11.4 &9.8 &9.7 &9.2 \\
skewness &-0.3 &1.2 &NA &1.4 &1.1 &1.1 &1.1 &1.2 \\
Kurtosis &-1.3 &1.7 &NA &0.7 &0.7 &-0.4 &-0.4 &-0.2 \\
median &8.7 &0.9 &NA &3.4 &10.3 &3.9 &2.3 &4.5 \\
\bottomrule
\end{tabular}
\end{threeparttable}
\end{table}

\pagebreak
\section{Equidistribution Measures}
 \subsection{Equidistribution on non-degenerate intervals}
A sequence $(s_1,s_2,s_3,...)$ of real numbers is said to be 
equidistributed on a non-degenerate interval $[a,b]$ if for every 
subinterval $[c,d]$ of $[a,b]$ we have

\begin{equation}
\lim _{n\to\infty}\frac{\left|\{\,s_{1},\dots,s_{n}\,\}\cap[c,d]\right|}{n}=\frac{d-c}{b-a}.    
\end{equation}\label{Eq:equinondegenerate}

(Here,the notation $|{s1,...,sn} \cap[c,d]|$ denotes the number of 
elements, out of the first $n$ elements of the sequence,that are between 
$c$ and $d$.) \cite{kuipers2012uniform}

To apply Eq.\ref{Eq:equinondegenerate}, the following steps where followed:
\begin{itemize}
    \item Normalize $\mathbf{Md}$
    \item Consider an interval with $size\in \{0.1,0.2, 0.25, 0.333\}$ 
    \item measure how many subintervals ($n$) of $\mathbf{Md}$ fall within a given interval, such that $size \leq |Md_{k+1}-Md_{k}|\leq 2\times size$
    \item With the previous measure, obtain the proportion of those subintervals, with respect to the length $N$ of the entire sequence, that is, $n/N$.
\end{itemize}
The results are summarized in the matrix of Table \ref{Table:ResultsNonDegenerate}

\begin{table}[hbt!]
\begin{threeparttable}
\caption{Measured vs expected fraction of subintervals of $\mathbf{Md}$}
\label{Table:ResultsNonDegenerate}
\begin{tabular}{cc}\toprule
 \text{measured n/N} & \text{expected} \\
 0.0909039 & 0.1 \\
 0.166657 & 0.2 \\
 0.199989 & 0.25 \\
 0.249986 & 0.333 \\
\bottomrule
\end{tabular}
\end{threeparttable}
\end{table}

\subsection{Well - distributed sequence}

A sequence $(s_1, s_2, s_3, \cdots)$ of real numbers is said to be \emph{well - 
 distributed} on $[a, b]$ if for any subinterval $[c, d]$ of $[a, b]$ we have

\begin{equation}\label{Eq:WellDistributed}
\lim_{n\to\infty}\frac{\left| \{\, s_{k + 1}, \dots, s_ {k + n}\, \}
\cap[c,d] \right|}{n} = \frac{d - c}{b - a}
\end{equation}

uniformly in $k$ \cite{kuipers2012uniform}.
To apply Eq.\ref{Eq:WellDistributed}, the following steps where followed:
\begin{itemize}
    \item Normalize $\mathbf{Md}$
    \item Consider an interval with $size\in \{0.1,0.2, 0.25, 0.333\}$ 
    \item measure how many subintervals ($n$) of $\mathbf{Md}$ fall within a given interval, such that $size \leq |Md_{k+j}-Md_{k}|\leq 2\times size$, where $j\in {1,2,3,4,5}$
    \item With the previous measure, obtain the proportion of those subintervals, with respect to the length $N$ of the entire sequence, that is, $n/N$, and for every $j$.
\end{itemize}
The results are summarized in Table \ref{Table:WellDistributedResults}

\begin{table}[hbt!]
\begin{threeparttable}
\caption{Well-distributedness properties of $\mathbf{Md}$: Measured vs expected fraction of subintervals of $\mathbf{Md}$, with varying spacing $k$}
\label{Table:WellDistributedResults}
\begin{tabular}{ccccccc}\toprule
\text{k=} & 1 & 2 & 3 & 4 & 5 & \text{Expected} \\
 \text{} & 0.0909039 & 0.0908987 & 0.0908935 & 0.0908884 & 0.0908832 & 0.1 \\
 \text{} & 0.166657 & 0.166648 & 0.166638 & 0.166629 & 0.166619 & 0.2 \\
 \text{} & 0.199989 & 0.199977 & 0.199966 & 0.199954 & 0.199943 & 0.25 \\
 \text{} & 0.249986 & 0.249972 & 0.249957 & 0.249943 & 0.249929 & 0.333 \\
\bottomrule
\end{tabular}
\end{threeparttable}
\end{table}

\subsection{Weyl' s criterion}
Weyl' s criterion \cite{kuipers2012uniform} states that the sequence an is equidistributed \
modulo 1 if and only if for all non - zero integers $\ell$
\begin{equation}
    \lim _ {{n\to\infty}} {\frac {1} {n}}\sum _ {{j = 
      1}}^{{n}} e^{{2\pi i\ell a_ {j}}} = 0.
\end{equation}

A quantitative form of Weyl' s criterion is given by the 
Erd\"os-Tur\'an inequality (\cite{ErdosTuran1948}).

Let $(s1, s2,s3 ...)\in \mathbb{R}$ be a sequence . 
   The Erd\"os-Tur\'an inequality applied to the measure
\begin{equation}
    \frac{1}{m} \#\{1\leq j\leq m\, | \, s_ {j}\, {\mathrm {mod}}\, 
     1\in S\}, \quad S\subset[0, 1),
\end{equation}

yields the following bound for the discrepancy:

\begin{equation}
    \begin{aligned} D (m) &\left(= \sup _ {{0\leq a\leq b\leq 1}} {\Big 
|} m^{{-1}} \#\{1\leq j\leq m\, | \, a\leq s_ {j}\, {\mathrm {mod}}\, 
           1\leq b\} - (b - a) {\Big |} \right) \\[
     8 pt] &\leq C\left({\frac {1} {n}} + {\frac {1} {m}}\sum _ {{k =
            1}}^{n} {\frac {1} {k}}\left| \sum _ {{j = 
          1}}^{m} e^{{2\pi is_ {j} k}} \right| \right) .
  \end{aligned}
\end{equation}\label{Eq:ErdosTuranWeyl}

This inequality holds for arbitrary natural numbers m, n, and gives a 
quantitative form of Weyl' s criterion for equidistribution. In Eq.\ref{Eq:ErdosTuranWeyl}, $m=|\mathbf{Md}|$, and $n$ is a natural number set up in this work to review if there is convergence.

\begin{figure}[htb]
\centering\includegraphics[width=10cm]{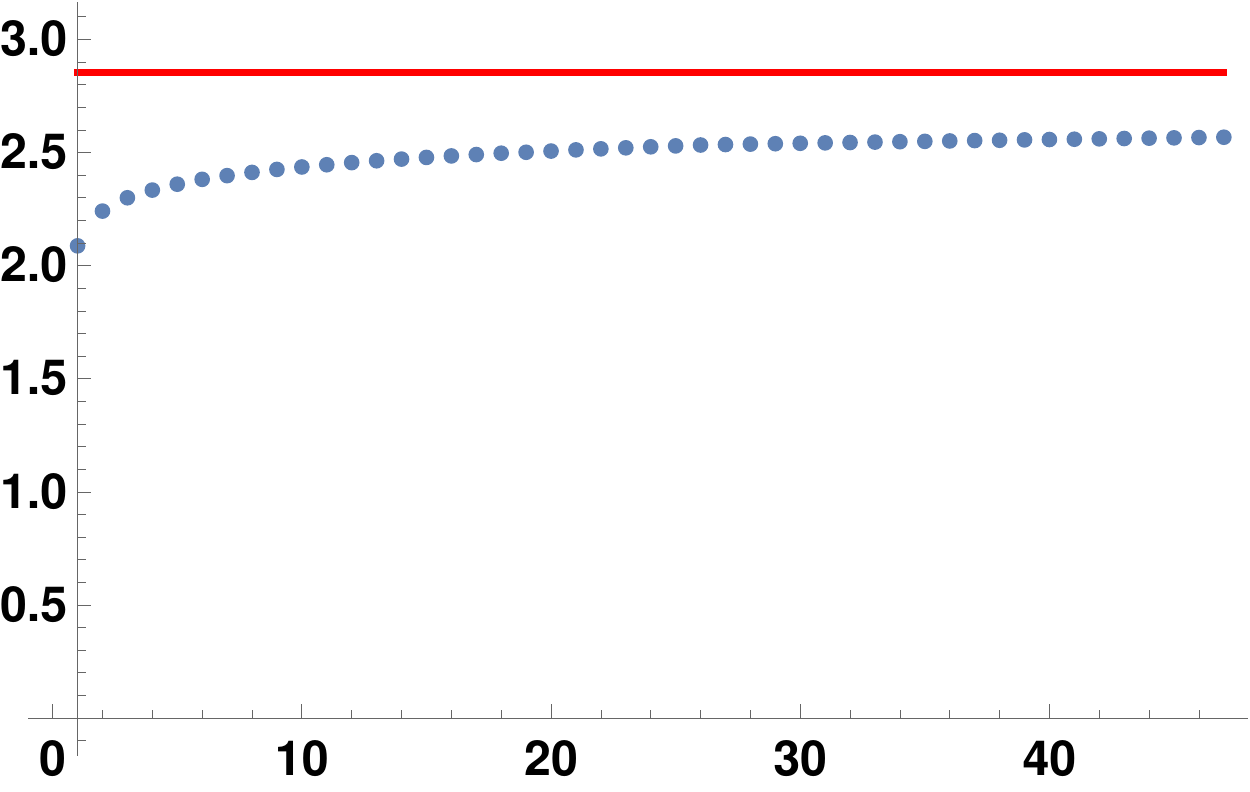}
\caption{Upper bound of $\mathbf{Md}\leq 2.85$, as a sequence of values as $m,n\to \infty$}\label{Figure:ErdosTuranWeyl}
\end{figure}

Applying Eq.\ref{Eq:ErdosTuranWeyl}, it was possible to determine that the upper bound of $\mathbf{Md}\leq 2.85$, as it is illustrated in Fig.\ref{Figure:ErdosTuranWeyl}

\begin{figure}[htb]
\centering\includegraphics[width=10cm]{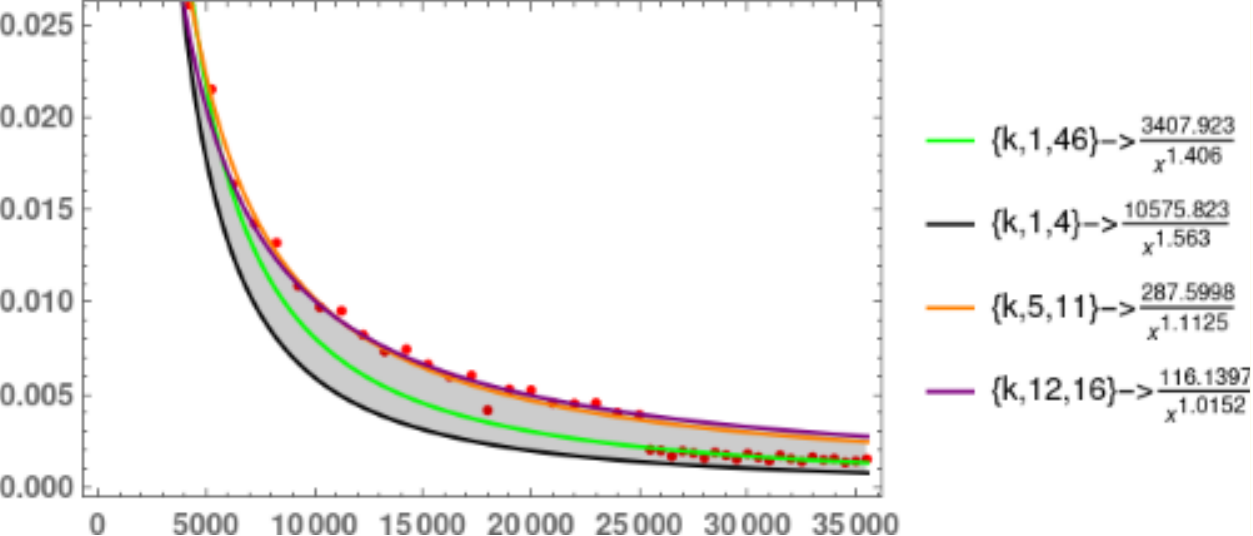}
\caption{Convergence of differences of upper bound for $\mathbf{Md}$, together with several fits of nonlinear models of the form $ax^b$ to illustrate these differences that the seem to be a Cauchy sequence}\label{Figure:ConvergenceCauchy}
\end{figure}

To verify that the upper bound shown in Fig.\ref{Figure:ErdosTuranWeyl} does converge, it was calculated the difference of the successor minus the predecessor of the sequence of partial sums/upper bounds. For the purposes of illustration, that set of differences was fit into several nonlinear models of the form $ax^b$, for several portions of the set, to show graphically that this set of differences are likely to be a Cauchy sequence, thus convergence is guaranteed. Both the set of differences and the non-linear fits are illustrated in Fig.\ref{Figure:ConvergenceCauchy} 

\begin{figure}[htb]
\centering\includegraphics[width=10cm]{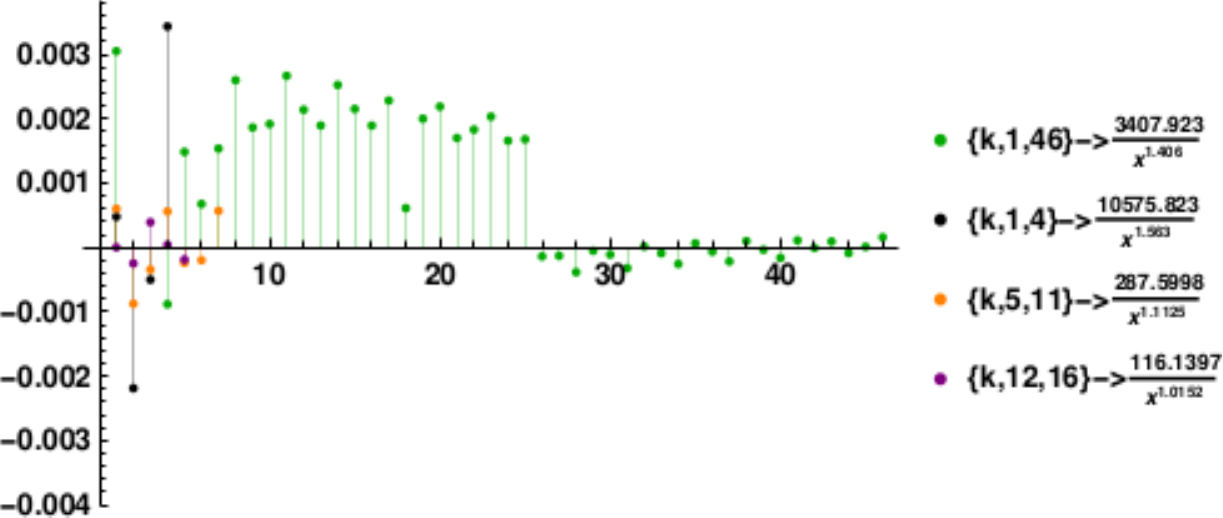}
\caption{Plot of residuals (the errors between the differences of partial sums and the non-linear fits of the form $ax^b$) showing the fullfill the criterium for being a Cauchy sequence, and thus the partial sums calculated using Eq.\ref{Eq:ErdosTuranWeyl} does converge.}\label{Figure:CauchyResiduals}
\end{figure}

To further assess whether or not the differences form a Cauchy sequence, it was calculated the set of residuals from the different non-linear fits. A Cauchy sequence is a sequence whose terms ultimately become arbitrarily close together, after sufficiently many initial terms 
have been discarded, and this is the Cauchy criterion for convergence of sequences. For each $\varepsilon >0$ there are only finitely many sequence members outside the epsilon 
tube. Regardless which $\varepsilon > 0$ we have, there is an index $N_0$, so that the sequence lies 
afterwards completely in the epsilon tube $(a - \varepsilon, a + \varepsilon)$ \cite{spivak1994calculus,courant1965introduction}. That the residuals make a Cauchy sequence is presented in Fig.\ref{Figure:CauchyResiduals}

\subsection{Conjecture related to the zeros of the Riemann's zeta function}
It has been proven (assuming the Riemann hypothesis) that the zeros of the Riemann's zeta function are equidistributed \cite{montgomery1973pair}. Since the zeroes of Riemann's zeta function are a way of encoding the position of prime numbers \cite{ZerosEncodePrimesWatKins}, it could be conjectured that, at least partially, this is the reason \emph{why} $\mathbf{Md}$ is an equidistributed sequence. 

For testing this conjecture, in \cite{katz1999zeroes} it is mentioned the the “pair correlation” of the zeroes, that they obey the laws of the Gaussian (orequivalently, Circular) Unitary Ensemble GUE. Specifically
\begin{equation}
    r_2 (GUE)(x)=1-\left(\frac{\sin \pi x}{\pi x}\right)^2
\end{equation}\label{Eq:GUE}

\begin{figure}[htb]
\centering\includegraphics[width=10cm]{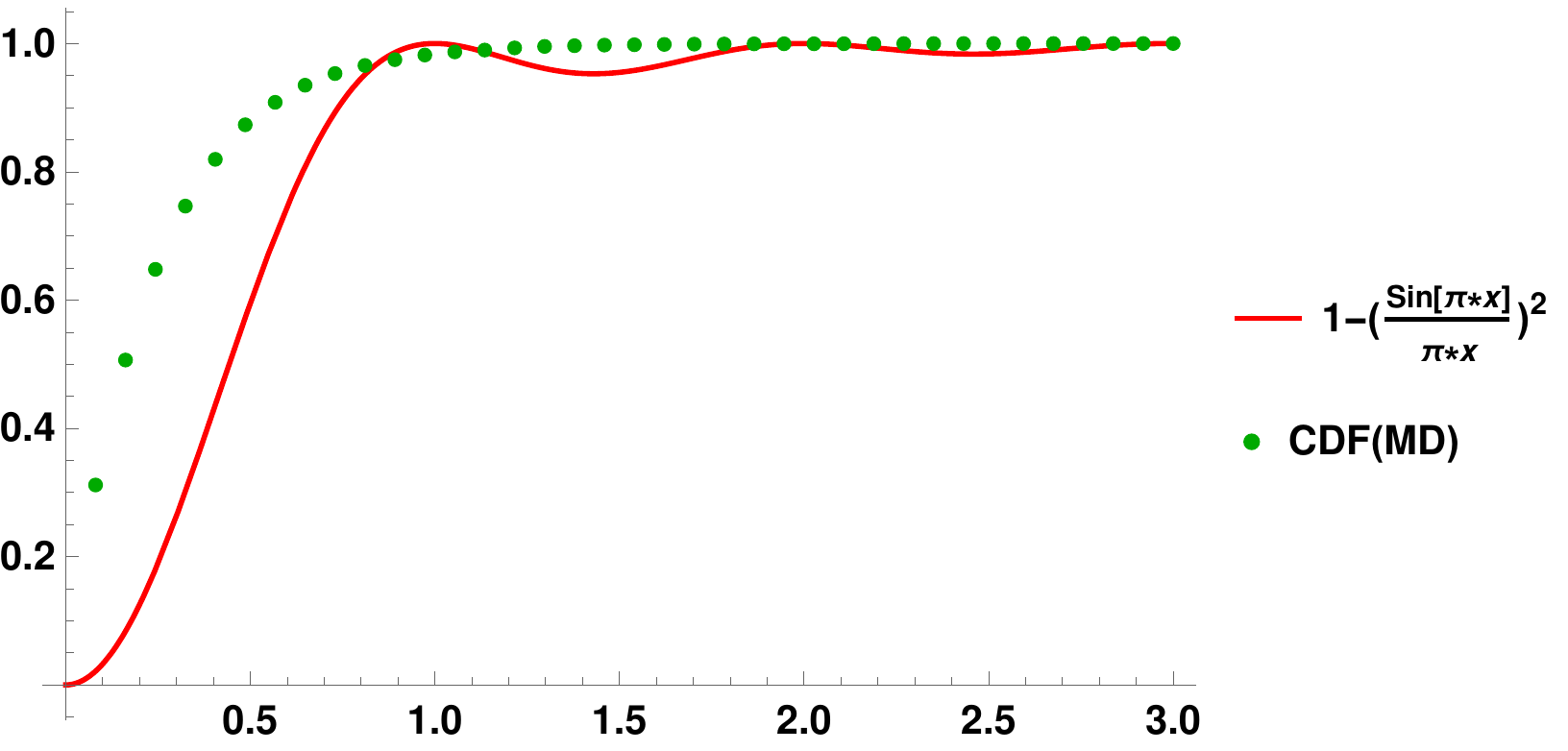}
\caption{Comparison of Eq.\ref{Eq:GUE} with the CFD of $\mathbf{MD}$}\label{Figure:GUECFDMD}
\end{figure}

\begin{figure}[htb]
\centering\includegraphics[width=10cm]{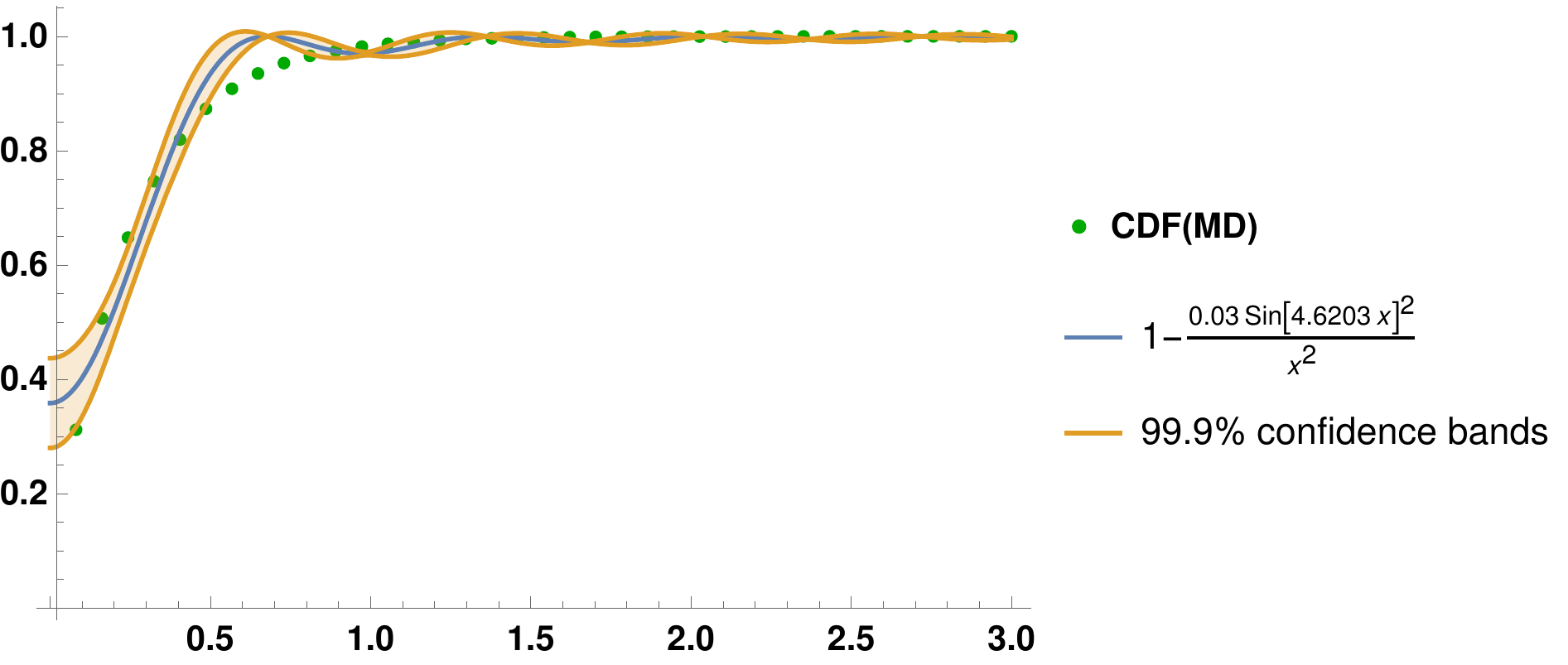}
\caption{Confidence bands at 99.9\% of the fitted parameters $a\approx 4.6203,b\approx 0.03$ to the model $1-\left(\frac{\sin \pi a x}{\pi b x}\right)^2$}\label{Figure:ErrorBands}
\end{figure}

\begin{figure}[htb]
\centering\includegraphics[width=10cm]{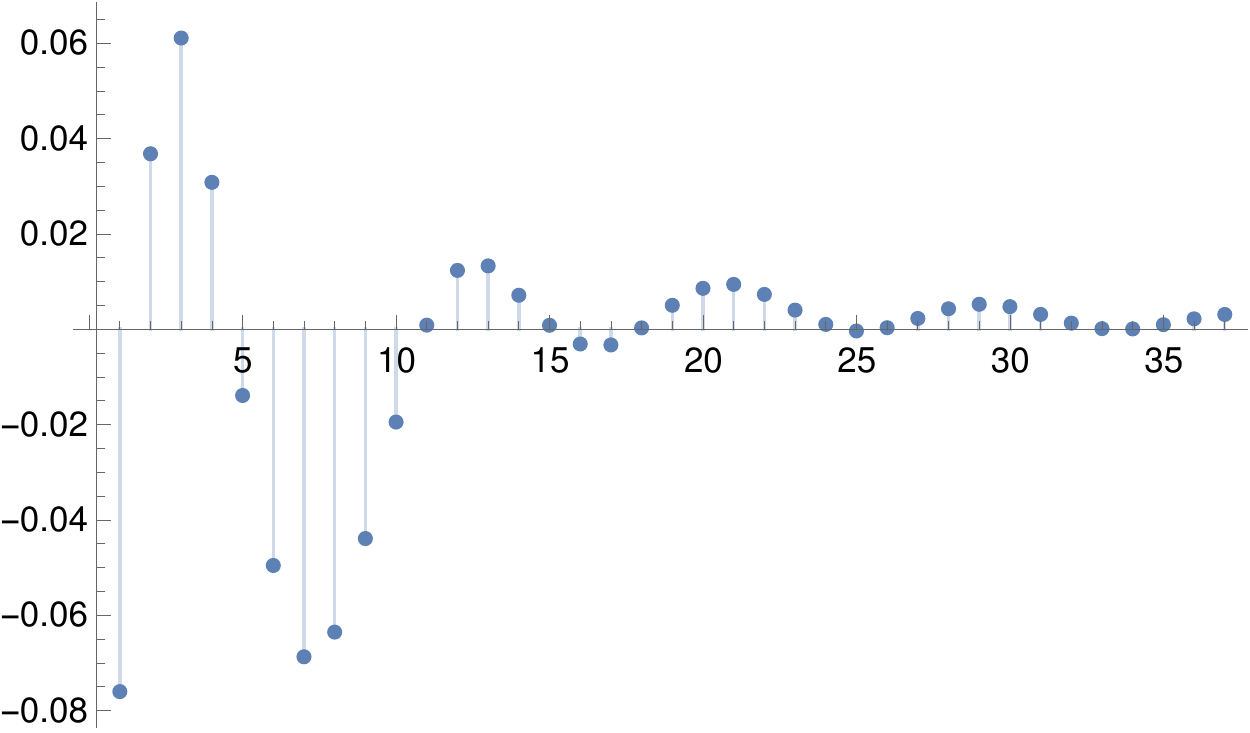}
\caption{Residuals of the model $1-\frac{0.03 \sin ^2(4.6203 x)}{x^2}$}\label{Figure:ResidualsNLM_CDF}
\end{figure}

Equation \ref{Eq:GUE} is compared to the Cumulative Distribution Function (CFD) of $\mathbf{MD}$ in Fig.\ref{Figure:GUECFDMD}; from this figure, it was thought that a non-linear model could be fitted, and the results are exposed in Figs.\ref{Figure:ErrorBands} and \ref{Figure:ResidualsNLM_CDF}. From these last figures, it could be conjectured that $\mathbf{Md}$ do seem to follow a model similar to Eq.\ref{Eq:GUE}, although not quite exactly.

\pagebreak
\section{Entropy}
In this work, the entropy of a random variable is considered to be the average level of "surprise", or "information" inherent to the variable's possible outcomes \cite{EntropyInfoTheoryWikipedia}. 

Given a discrete random variable $X$, which takes values in the alphabet ${\mathcal {X}}$ and is distributed according to $p: \mathcal{X}\to[0, 1]$: 
\begin{equation}
    \Eta(X) := -\sum_{x \in \mathcal{X}} p(x) \log p(x) = \mathbb{E}[-\log p(X)] ,
\end{equation}
It is conjectured that Low-Discrepancy Sequences contain "just enough" information for integration, as opposite of a numerical scheme that requires partitioning the interval at regular intervals (less "surprise") or a Monte Carlo scheme, that requires many functional evaluations (more "surprise"). In Table \ref{Table:Entropies} the values of entropies for several numerical system is presented. Notice that $\mathbf{Md}$ occupy a relative intermediate position concerning the numerical value of its Shannon entropy.

\begin{table}[!htp]\centering
\caption{Entropies of selected sequences}\label{Table:Entropies}
\scriptsize
\begin{tabular}{p{2in}p{2in}}\toprule
H(x) (entropy, or "degree of surprise") &Set measured \\\midrule
0 &(no surprise, like a crystal, no choice) The first million differences between random integers in the range 1 to 0ne million \\
0 &(no surprise, like a crystal, no choice) The first million differences between random reals in the range 1 to 0ne million \\
0.30103 &(This is the flipping of a coin one million times) \\
0.30103 &The first million digits of Pi \\
0.30103 &The first million digits of E \\
0.30103 &the first million digits of Sqrt[71] \\
0.778151 &(rolling of one dice) \\
1.26895 &The difference between the first million primes \\
1.63164 &The difference of exponents of the first 47 Mersenne Primes \\
1.6721 &The exponent of the first 47 Mersenne Primes \\
2.13791 &Metadistances {6,6} \\
2.89478 &(Entropy of "fiery" poem (\cite{blake1793america}), 20 words) \\
2.99573 &(entropy of a random choice of 20 common words) \\
3.00911 &(entropy of poem "no man is an Island" (\cite{donne1623whom}), 79 words) \\
4.36945 &(entropy of random choice of common words, 79 words) \\
5.75118 &(entropy of a million random integers in the range 1 to one million) \\
5.86415 &(entropy of a million random integers in the range minus one million to one million) \\
6 &(entropy of a million random reals in the range minus one million to one million) \\
6 &the first million primes \\
\bottomrule
\end{tabular}
\end{table}
This recent numerical findings and its relationship to $\mathbf{Md}$ are supported by previous analytic work \cite{tapia2016studies}.

\pagebreak

\section{Conclusion}
A revision of previous work was performed, to find new evidence that $\mathbf{Md}$ are an Equidistributed Sequence. It was conjectured that Ramanujan's master theorem serves as better argument for explaining the convergence of the solution while making numerical quadratures. The new statistics about previous results, show that $\mathbf{Md}$ has the most symmetrical distribution of errors, and the least amount of extreme outliers for obtaining a result.

The $\mathbf{Md}$ comply with all the equidistribution measures here studied, thus confirming this set as a Low-Discrepancy Sequence and an Equidistributed Sequence; the reason \emph{why} they are equidistributed it is conjectured may be related to the equidistribution of the zeros of Riemann's zeta function, and some numerical evidence is presented to show statistical similarities with the GUE model.

Finally, it is show a comparison of the entropy for several numerical systems or encoding, whereupon $\mathbf{Md}$ results having an intermediate level of entropy, possibly suggesting that $\mathbf{Md}$ have the "right" amount of information for numerical quadrature, as opposite to the amount of information required for Monte Carlo or other numerical methods.


\paragraph{Competing Interests}

None.


\printbibliography

\end{document}